\documentclass[11pt,leqno]{article}
\usepackage[mathscr]{eucal}
\usepackage{amsmath}
\usepackage{amssymb}
\usepackage{amsfonts}
\usepackage{latexsym}
\usepackage{theorem}
\theoremstyle{plain}

\newtheorem{Theorem}{Theorem}[section]

\newtheorem{Lemma}{Lemma}[section]
\newtheorem{Proposition}{Proposition}[section]

\theorembodyfont{\rmfamily}
\newtheorem{Remark}{Remark}[section]
\newtheorem{Example}{Example}[section]
\newtheorem{Definition}{Definition}[section]

\numberwithin{equation}{section}

\title{Minimal surfaces in a certain
$3$-dimensional homogeneous spacetime}
\author{Sungwook Lee}
\date{}
\begin{document}
\maketitle

\begin{abstract}
 The 2-parameter family of certain homogeneous Lorentzian\\ 3-manifolds which includes Minkowski 3-space, de Sitter 3-space, and Minkowski motion group is considered. Each homogeneous Lorentzian 3-manifold in the 2-parameter family has a solvable Lie group structure with left invariant metric. A generalized integral representation formula which is the unification of representation formulas for minimal timelike surfaces in those homogeneous Lorentzian 3-manifolds is obtained.  The normal Gau{\ss} map of minimal timelike surfaces in those homogeneous Lorentzian 3-manifolds and its harmonicity are discussed.
\end{abstract}


{\bf 2000 MR Subject Classification:} 53A10, 53C30, 53C42, 53C50
\section*{Introduction}
In \cite{Lee}, the author considered the $2$-parameter family of 3-dimensional homogeneous spacetimes $({\mathbb R}^3(x^0,x^1,x^2),g_{(\mu_1,\mu_2)})$ with Lorentzian metric
$$g_{(\mu_1,\mu_2)}=-(dx^{0})^{2}+ e^{-2\mu_1x^0}(dx^{1})^{2}
+e^{-2\mu_2x^0}(dx^{2})^{2}.$$
Every homogeneous Lorentzian manifold in this family can be represented as a solvable matrix Lie
group with left invariant metric
$$
G(\mu_1,\mu_2)=
\left\{\left(
\begin{array}{cccc}
1 & 0 & 0 & x^{0}\\
0 & e^{\mu_1x^0} & 0 & x^{1}\\
0 & 0 & e^{\mu_2x^0} & x^{2}\\
0 & 0 & 0 & 1
\end{array}
\right)
\
\Biggr
\vert
\
x^0, x^1,x^2
\in \mathbb{R}
\right
\}.
$$
As special cases, this family of homogeneous Lorentzian 3-manifolds include Minkowski $3$-space $\mathbb{E}^3_1$, de Sitter $3$-space $\mathbb{S}^3_1(c^2)$ of constant sectional curvature $c^2$ as a warped product model, and $\mathbb{S}^2_1(c^2)\times\mathbb{E}^1$, the direct product of de Sitter 2-space $\mathbb{S}^2_1(c^2)$ of constant curvature $c^2$ and the real line $\mathbb{E}^1$. (In fact, Minkowski $3$-space and de Sitter $3$-space are the only homogeneous Lorentzian 3-manifolds in this family that have a constant sectional curvature.) These three spaces may be considered as Lorentzian counterparts of Euclidean 3-space $\mathbb{E}^3$, hyperbolic 3-space $\mathbb{H}^3(-c^2)$ of constant sectional curvature $-c^2$, and $\mathbb{H}^2(-c^2)\times\mathbb{E}^1$, the direct product of hyperbolic plane $\mathbb{H}^2(-c^2)$ of constant curvature $-c^2$ and the real line $\mathbb{E}^1$, respectively, of Thurston's eight model geometries \cite{Thurston}. In \cite{Lee}, the author obtained a generalized integral representation formula which is the unification of representation formulas for maximal spacelike surfaces in those homogeneous Lorentzian 3-manifolds. In particular, the generalized integral formula includes Weierstra{\ss} representation formula for maximal spacelike surfaces in Minkowski 3-space studied independently by O. Kobayashi \cite{Kobayashi} and L. McNertney \cite{Mc}, and Weierstra{\ss} representation formula for maximal spacelike surfaces in de Sitter 3-space.

In this paper, the author obtains a generalized integral representation formula which is the unification of representation formulas for minimal timelike surfaces in those homogeneous Lorentzian 3-manifolds. In particular, the generalized integral formula includes Weierstra{\ss} representation formula for minimal timelike surfaces in Minkowski 3-space (\cite{Inoguchi-Toda}, \cite{Lee2}) and Weierstra{\ss} representation formula for minimal timelike surfaces in de Sitter 3-space. The harmonicity of the normal Gau{\ss} map of minimal timelike surfaces in $G(\mu_1,\mu_2)$ is also discussed. It is shown that Minkowski $3$-space $G(0,0)$, de Sitter $3$-space $G(c,c)$, and Minkowski motion group $G(c,-c)$ are the only homogeneous Lorentzian $3$-manifolds among the 2-parameter family members $G(\mu_1,\mu_2)$ in which the (projected) normal Gau{\ss} map of minimal timelike surfaces is harmonic. The harmonic map equations for those cases are obtained.
\section{Solvable Lie group}
\label{sec:sol}
In this section, we study the two-parameter family of certain homogeneous
Lorentzian $3$-manifolds.

Let us consider the two-parameter family of homogeneous Lorentzian $3$-manifolds
\begin{equation}\label{1.1}
\left\{(\mathbb{R}^{3}(x^{0},x^{1},x^{2}), g_{(\mu_1,\mu_2)})\
\vert \
(\mu_1,\mu_2)\in \mathbb{R}^{2}
\right \},
\end{equation}
where the metric $g_{(\mu_1,\mu_2)}$ is
defined by
\begin{equation}
g_{(\mu_1,\mu_2)}:=-(dx^{0})^{2}+
e^{-2\mu_1x^0}(dx^{1})^{2}
+e^{-2\mu_2x^0}(dx^{2})^{2}.
\end{equation}
\begin{Proposition}
Each homogeneous space $(\mathbb{R}^3,g_{(\mu_1,\mu_2)})$ is
isometric to the following solvable matrix Lie group:
$$
G(\mu_1,\mu_2)=
\left\{\left(
\begin{array}{cccc}
1 & 0 & 0 & x^{0}\\
0 & e^{\mu_1x^0} & 0 & x^{1}\\
0 & 0 & e^{\mu_2x^0} & x^{2}\\
0 & 0 & 0 & 1
\end{array}
\right)
\
\Biggr
\vert
\
x^0, x^1,x^2
\in \mathbb{R}
\right
\}
$$
with left invariant metric. The group operation on $G(\mu_1,\mu_2)$ is the ordinary matrix multiplication and the corresponding group operation on $(\mathbb{R}^3,g_{(\mu_1,\mu_2)})$ is given by
$$
(x^0,x^1,x^2)\cdot (\tilde{x}^0,\tilde{x}^1,\tilde{x}^2)=
(x^0+\tilde{x}^0,x^1+e^{\mu_1x^0}\tilde{x}^1,
x^2+e^{\mu_2x^0}\tilde{x}^2).
$$
\end{Proposition}
{\it Proof.} For $\tilde{a}=(a^0,a^1,a^2)\in G(\mu_1,\mu_2)$, denote by $L_{\tilde{a}}$ the left translation by $\tilde{a}$. Then 
$$L_{\tilde{a}}(x^0,x^1,x^2)=(x^0+a^0,e^{\mu_1a^0}x^1+a^1,e^{\mu_2a^0}x^2+a^2)$$
and
\begin{align*}
L^\ast_{\tilde{a}}g_{(\mu_1,\mu_2)}=&-\{d(x^0+a^0)\}^2+e^{-2\mu_1(x^0+a^0)}\{d(e^{\mu_1a^0}x^1+a^1)\}^2\\
&+e^{-2\mu_2(x^0+a^0)}\{d(e^{\mu_2a^0}x^2+a^2)\}^2\\
=&-(dx^0)^2+e^{-2\mu_1x^0}(dx^1)^2+e^{-2\mu_2x^0}(dx^2)^2.
\end{align*}
Q.E.D.

The Lie algebra $\mathfrak{g}(\mu_1,\mu_2)$ is
given by
\begin{equation}
\mathfrak{g}(\mu_1,\mu_2)=
\left\{\left(
\begin{array}{cccc}
0 & 0 & 0 & y^{0}\\
0 & \mu_1y^0 & 0 & y^{1}\\
0 & 0 & \mu_2y^0 & y^{2}\\
0 & 0 & 0 & 0
\end{array}
\right)
\
\Biggr
\vert
\
y^0,y^1,y^2
\in \mathbb{R}
\right
\}.
\end{equation}
We take the following basis $\{E_0,E_1,E_2\}$
of $\mathfrak{g}(\mu_1,\mu_2)$:
\begin{equation}
\label{eq:basis}
E_{0}= \left(
\begin{array}{cccc}
0 & 0 & 0 & 1\\
0 & \mu_1 & 0 & 0\\
0 & 0 & \mu_2 & 0\\
0 & 0 & 0 & 0
\end{array}
\right), E_{1}= \left(
\begin{array}{cccc}
0 & 0 & 0 & 0\\
0 & 0 & 0 & 1\\
0 & 0 & 0 & 0\\
0 & 0 & 0 & 0
\end{array}
\right),
E_{2}=
\left(
\begin{array}{cccc}
0 & 0 & 0 & 0\\
0 & 0 & 0 & 0\\
0 & 0 & 0 & 1\\
0 & 0 & 0 & 0
\end{array}
\right).
\end{equation}
Then the commutation relation of $\mathfrak{g}(\mu_1,\mu_2)$ is
given by
$$
[E_1,E_2]=0,\
[E_2,E_0]=-\mu_{2}E_{2},\
[E_0,E_1]=\mu_{1}E_{1}.
$$
The left translation of $E_0,E_1,E_2$ are the vector fields $e_0=\frac{\partial}{\partial x^0}$, $e_1=e^{\mu_1x^0}\frac{\partial}{\partial x^1}$, $e_2=e^{\mu_2x^0}\frac{\partial}{\partial x^2}$, respectively such that
\begin{align*}
\langle e_0,e_0\rangle&=-1,\ \langle e_1,e_1\rangle=\langle e_2,e_2\rangle=1,\\
\langle e_i,e_j\rangle&=0\ \mbox{if}\ i\ne j.
\end{align*}
That is, $\{e_0,e_1,e_2\}$ forms a Lorentzian frame field on $(\mathbb{R}^3,g_{(\mu_1,\mu_2)})$. Hence we see that $\{E_0,E_1,E_2\}$ forms an orthonormal basis for $\mathfrak{g}(\mu_1,\mu_2)$.

For $X\in{\mathfrak g}(\mu_1,\mu_2)$, denote by ${\rm ad}(X)^*$ the
\emph{adjoint} operator of ${\rm ad}(X)$ i.e. it is defined by
the equation
$$\langle{\rm ad}(X)(Y),Z\rangle=\langle Y,{\rm ad}(X)^*(Z)\rangle$$
for any $Y,Z\in{\mathfrak g}(\mu_1,\mu_2)$. Here ${\rm ad}(X)(Y)=[X,Y]$ for $X,Y\in\mathfrak{g}$. Let $U$ be the symmetric
bilinear operator on ${\mathfrak g}(\mu_1,\mu_2)$ defined by
$$U(X,Y):=\frac{1}{2}\{{\rm ad}(X)^*(Y)+{\rm ad}(Y)^*(X)\}.$$
\begin{Lemma}
Let $\{E_0,E_1,E_2\}$ be the orthonormal basis for ${\mathfrak
g}(\mu_1,\mu_2)$ defined in \eqref{eq:basis}. Then
\begin{align*}
U(E_0,E_0)&=0,\ U(E_1,E_1)=\mu_1E_0,\ U(E_2,E_2)=\mu_2E_0,\\
U(E_0,E_1)&=\frac{\mu_1}{2}E_1,\ U(E_1,E_2)=0,\
U(E_2,E_0)=\frac{\mu_2}{2}E_2.
\end{align*}
\end{Lemma}

Let $\mathfrak{D}$ be a simply connected domain and $\varphi:\mathfrak{D}\longrightarrow G(\mu_1,\mu_2)$ an immersion. $\varphi$ is said to be \emph{timelike} if the induced metric $I$ on $\mathfrak{D}$ is Lorentzian. The induced Lorentzian metric $I$ determines a Lorentz conformal structure $\mathcal{C}_I$ on $\mathfrak{D}$. Let $(t,x)$ be a Lorentz isothermal coordinate system with respect to the conformal structure $\mathcal{C}_I$. Then the first fundamental form $I$ is written in terms of $(t,x)$ as
\begin{equation}
I=e^\omega (-dt^2+dx^2).
\end{equation}
The conformality condition is given in terms of $(t,x)$ by
\begin{equation}
\label{eq:conf}
\begin{aligned}
\left\langle\frac{\partial\varphi}{\partial t},\frac{\partial\varphi}{\partial x}\right\rangle&=0,\\
-\left\langle\frac{\partial\varphi}{\partial t},\frac{\partial\varphi}{\partial t}\right\rangle&=\left\langle\frac{\partial\varphi}{\partial x},\frac{\partial\varphi}{\partial x}\right\rangle=e^\omega.
\end{aligned}
\end{equation}
A conformal timelike surface is called a \emph{Lorentz surface}. Let $u:=t+x$ and $v:=-t+x$. Then $(u,v)$ defines a null coordinate system with respect to the conformal structure $\mathcal{C}_I$. The first fundamental form $I$ is written in terms of $(u,v)$ as
\begin{equation}
I=e^\omega dudv.
\end{equation}
The partial derivatives $\frac{\partial\varphi}{\partial u}$ and $\frac{\partial\varphi}{\partial v}$ are computed to be
\begin{equation}
\frac{\partial\varphi}{\partial u}=\frac{1}{2}\left(\frac{\partial\varphi}{\partial t}+\frac{\partial\varphi}{\partial x}\right),\ \frac{\partial\varphi}{\partial v}=\frac{1}{2}\left(-\frac{\partial\varphi}{\partial t}+\frac{\partial\varphi}{\partial x}\right).
\end{equation}
The conformality condition \eqref{eq:conf} can be written in terms of null coordinates as
\begin{equation}
\begin{aligned}
\left\langle\frac{\partial\varphi}{\partial u},\frac{\partial\varphi}{\partial u}\right\rangle&=\left\langle\frac{\partial \varphi}{\partial v},\frac{\partial\varphi}{\partial v}\right\rangle=0,\\
\left\langle\frac{\partial\varphi}{\partial u},\frac{\partial\varphi}{\partial v}\right\rangle&=\frac{1}{2}e^\omega.
\end{aligned}
\end{equation}
\begin{Definition}
Let ${\mathfrak D}(t,x)$ be a simply connected domain. A smooth timelike immersion $\varphi: {\mathfrak D}\longrightarrow G(\mu_1,\mu_2)$ is said to be \emph{harmonic} if it is a critical point of the energy functional\footnote{This is an analogue of the Dirichlet energy.}
\begin{equation}
E(\varphi)=\int_{\mathfrak{D}}e(\varphi)dtdx,
\end{equation} 
where $e(\varphi)$ is the \emph{energy density} of $\varphi$
\begin{equation}
e(\varphi)=\frac{1}{2}\left\{-\left|\varphi^{-1}\frac{\partial\varphi}{\partial t}\right|^2+\left|\varphi^{-1}\frac{\partial\varphi}{\partial x}\right|^2\right\}.
\end{equation}
$\left|\varphi^{-1}\frac{\partial\varphi}{\partial t}\right|^2=\left\langle\varphi^{-1}\frac{\partial\varphi}{\partial t},\varphi^{-1}\frac{\partial\varphi}{\partial t}\right\rangle<0$ and $\left|\varphi^{-1}\frac{\partial\varphi}{\partial x}\right|^2=\left\langle\varphi^{-1}\frac{\partial\varphi}{\partial x},\varphi^{-1}\frac{\partial\varphi}{\partial x}\right\rangle>0$, so $e(\varphi)>0$ and hence $E(\varphi)\geq 0$.
\end{Definition}
\begin{Lemma}

Let ${\mathfrak D}$ be a simply connected domain. A smooth timelike immersion
$\varphi: {\mathfrak D}\longrightarrow G(\mu_1,\mu_2)$ is harmonic
if and only if it satisfies the wave equation
\begin{equation}
\label{eq:harm}
\begin{aligned}
-&\frac{\partial}{\partial t}\left(\varphi^{-1}\frac{\partial\varphi}{\partial t}\right)+\frac{\partial}{\partial x}\left(\varphi^{-1}\frac{\partial\varphi}{\partial x}\right)\\
&-\left\{-\mathrm{ad}\left(\varphi^{-1}\frac{\partial\varphi}{\partial t}\right)^\ast\left(\varphi^{-1}\frac{\partial\varphi}{\partial t}\right)+\mathrm{ad}\left(\varphi^{-1}\frac{\partial\varphi}{\partial x}\right)^\ast\left(\varphi^{-1}\frac{\partial\varphi}{\partial x}\right)\right\}=0.
\end{aligned}
\end{equation}
\end{Lemma}
{\it Proof.} Let $\varphi_s$, $s\in (-\epsilon,\epsilon)$ be a smooth variation of $\varphi=\varphi_0$ such that $\varphi_s|_{\partial\mathfrak{D}}=\varphi|_{\partial\mathfrak{D}}$, where $\partial\mathfrak{D}$ is the boundary $\mathfrak{D}$. Let
$$\Lambda=\frac{d}{ds}(\varphi^{-1}\varphi_s)\vert_{s=0}:\mathfrak{D}\longrightarrow\mathfrak{g}(\mu_1,\mu_2).$$
\begin{align*}
\frac{d}{ds}E(\varphi_s)\vert_{s=0}&=\int_{\mathfrak{D}}\left\{-\left\langle\frac{d}{ds}\left(\varphi^{-1}\frac{\partial\varphi}{\partial t}\right)\vert_{s=0},\varphi^{-1}\frac{\partial\varphi}{\partial t}\right\rangle\right.\\
&+\left.\left\langle\frac{d}{ds}\left(\varphi^{-1}\frac{\partial\varphi}{\partial x}\right)\vert_{s=0},\varphi^{-1}\frac{\partial\varphi}{\partial x}\right\rangle\right\}dtdx\\
&=\int_{\mathfrak{D}}
\left\{-\left\langle\left[\varphi^{-1}\frac{\partial\varphi}{\partial t},\Lambda\right]+\frac{\partial\Lambda}{\partial t},\varphi^{-1}\frac{\partial\varphi}{\partial t}\right\rangle\right.\\
&+\left.\left\langle\left[\varphi^{-1}\frac{\partial\varphi}{\partial x},\Lambda\right]+\frac{\partial\Lambda}{\partial x},\varphi^{-1}\frac{\partial\varphi}{\partial x}\right\rangle\right\}dtdx\\
&=\int_{\mathfrak{D}}\left\langle\Lambda,-\frac{\partial}{\partial t}\left(\varphi^{-1}\frac{\partial\varphi}{\partial t}\right)+\frac{\partial}{\partial x}\left(\varphi^{-1}\frac{\partial\varphi}{\partial x}\right)\right.\\
&+\left.\mathrm{ad}\left(\varphi^{-1}\frac{\partial\varphi}{\partial t}\right)^\ast\left(\varphi^{-1}\frac{\partial\varphi}{\partial t}\right)-\mathrm{ad}\left(\varphi^{-1}\frac{\partial\varphi}{\partial x}\right)^\ast\left(\varphi^{-1}\frac{\partial\varphi}{\partial x}\right)\right\rangle\\
&dtdx.
\end{align*}
This completes the proof.

In terms of null coordinates $u,v$, the
wave equation \eqref{eq:harm} can be written as
\begin{equation}
\label{eq:harm2}
\frac{\partial}{\partial u}\left(\varphi^{-1}\frac{\partial\varphi}{\partial
  v}\right)+\frac{\partial}{\partial
  v}\left(\varphi^{-1}\frac{\partial\varphi}{\partial u}\right)-2U\left(\varphi^{-1}\frac{\partial\varphi}{\partial
  u},\varphi^{-1}\frac{\partial\varphi}{\partial v}\right)=0.
\end{equation}
Let $\varphi^{-1}d\varphi=\alpha'du+\alpha''dv$. Then the equation
\eqref{eq:harm2} is equivalent to
\begin{equation}
\label{eq:harm3}
\alpha'_v+\alpha''_u=2U(\alpha',\alpha'').
\end{equation}
The Maurer-Cartan equation is given by
\begin{equation}
\label{eq:m-c}
\alpha'_v-\alpha''_u=[\alpha',\alpha''].
\end{equation}
The equations \eqref{eq:harm3} and \eqref{eq:m-c} can be combined to a single
equation
\begin{equation}
\label{eq:int}
\alpha'_v=U(\alpha',\alpha'')+\frac{1}{2}[\alpha',\alpha''].
\end{equation}
The equation \eqref{eq:int} is both the integrability condition for the
differential equation $\varphi^{-1}d\varphi=\alpha'du+\alpha''dv$ and the
condition for $\varphi$ to be a harmonic map.

The Levi-Civita connection $\nabla$ of
$G(\mu_1,\mu_2)$ is computed to be
\begin{align*}
\nabla_{e_0}e_{0}&=0,\ \nabla_{e_0}e_{1}=-\mu_1e_1,\ \nabla_{e_0}e_{2}=-\mu_2e_2,\\
\nabla_{e_1}e_{0}&=-\mu_1e_1,\
\nabla_{e_1}e_{1}=-\mu_1e_0,\ \nabla_{e_1}e_{2}=0,\\
\nabla_{e_2}e_{0}&=-\mu_2e_2,\ \nabla_{e_2}e_{1}=0,\
\nabla_{e_2}e_{2}=-\mu_2e_0.
\end{align*}
Let $K(e_i,e_j)$ denote the sectional curvature of $G(\mu_1,\mu_2)$
with respect to the tangent plane spanned by $e_i$ and $e_j$ for
$i,j=0,1,2$. Then
\begin{equation}
\label{eq:curv}
\begin{aligned}
K(e_0,e_1)&=g^{00}R^1_{010}=\mu_1^2,\\
K(e_1,e_2)&=g^{11}R^2_{121}=\mu_1\mu_2,\\
K(e_0,e_3)&=g^{00}R^3_{030}=\mu_2^2,
\end{aligned}
\end{equation}
where $g_{ij}=g_{(\mu_1,\mu_2)}(e_i,e_j)$ denotes the metric tensor
of $G(\mu_1,\mu_2)$. Hence, we see that $G(\mu_1,\mu_2)$ has
constant sectional curvature if and only if
$\mu_1^2=\mu_2^2=\mu_1\mu_2$. If $c:=\mu_1=\mu_2$, then
$G(\mu_1,\mu_2)$ is locally isometric to ${\mathbb S}^3_1(c^2),$ the
de Sitter $3$-space of constant sectional curvature $c^2$. (See
Example \ref{ex:desitter} and Remark \ref{rem:desitter}.) If
$\mu_1=-\mu_2$, then $\mu_1=\mu_2=0$, so
$G(\mu_1,\mu_2)=G(0,0)$ is locally isometric to ${\mathbb E}^3_1$ (Example
\ref{ex:minkowski}).
\begin{Example}(Minkowski $3$-space)
\label{ex:minkowski}
The Lie group $G(0,0)$ is isomorphic and
isometric to the Minkowski $3$-space
$$
\mathbb{E}^3_{1}=(\mathbb{R}^3(x^0,x^1,x^2),+)
$$
with the metric $-(dx^0)^2+(dx^1)^2+(dx^2)^2$.
\end{Example}

\begin{Example}(de Sitter $3$-space)
\label{ex:desitter}
Take $\mu_1=\mu_2=c \not=0$. Then $G(c,c)$ is
the flat chart model of the de Sitter $3$-space:
$$
\mathbb{S}^{3}_{1}(c^2)_{+}=(\mathbb{R}^3(x^0,x^1,x^2),
-(dx^0)^{2}+e^{-2cx^0}\{(dx^{1})^{2}+(dx^{2})^{2}\}).
$$
\end{Example}

\begin{Remark}
\label{rem:desitter} Let $\mathbb{E}^{4}_1$ be the Minkowski
$4$-space. The natural Lorentzian metric $\langle \cdot,\cdot
\rangle$ of $\mathbb{E}^{4}_1$ is expressed as
$$
\langle \cdot,\cdot\rangle=-(du^{0})^{2}+(du^{1})^{2}+
(du^2)^{2}+(du^3)^2.
$$
in terms of natural coordinate system
$(u^0,u^1,u^2,u^3)$. The {\it de Sitter $3$-space} $\mathbb{S}^3_1(c^2)$ of constant
sectional curvature $c^2>0$ is realized as the hyperquadric in
$\mathbb{E}^{4}_{1}$:
$$
\mathbb{S}^3_1(c^2)= \{(u^0,u^1,u^2,u^3)\in\mathbb{E}^{4}_1: \
-(u^{0})^{2}+(u^{1})^{2}+(u^2)^2+(u^3)^2=1/c^2\}.
$$
The de Sitter $3$-space $\mathbb{S}^{3}_1(c^2)$ is divided into the
following three regions:
\begin{align*}
\mathbb{S}^{3}_{1}(c^2)_{+}&= \{(u^0,u^1,u^2,u^3)\in
\mathbb{S}^{3}_1(c^2): \ c(u^{0}+u^{1})>0\};\\
\mathbb{S}^{3}_{1}(c^2)_{0}&=
\{(u^0,u^1,u^2,u^3)\in\mathbb{S}^{3}_1(c^2): \ u^{0}+u^{1}=0\};\\
\mathbb{S}^{3}_{1}(c^2)_{-}&=
\{(u^0,u^1,u^2,u^3)\in\mathbb{S}^{3}_1(c^2): \ c(u^{0}+u^{1})<0\}.
\end{align*}
$\mathbb{S}^3_1(c^2)$ is the disjoint union
$\mathbb{S}^3_1(c^2)_{+}\dotplus\mathbb{S}^3_1(c^2)_{0}\dotplus
\mathbb{S}^3_1(c^2)_{-}$ and $\mathbb{S}^{3}_{1}(c^2)_{\pm}$ are
diffeomorphic to $\mathbb{R}^3$. Let us introduce a local coordinate system $(x^0,x^1,x^2)$ by
$$
x^0=\frac{1}{c}\log c(u^0+u^1),\ \
x^{j}=\frac{u^{j+1}}{c(u^0+u^1)},\ (j=1,2).
$$
This local coordinate system is defined on
$\mathbb{S}^3_1(c^2)_{+}$. The induced metric of
$\mathbb{S}^3_1(c^2)_+$ is expressed as:
$$
g_c:=-(dx^0)^2+e^{2cx^0}\{(dx^1)^2+(dx^2)^2\}.
$$
The chart $(\mathbb{S}^3_1(c^2)_{+},g_c)$ is traditionally called
the {\it flat chart} of $\mathbb{S}^3_1(c^2)$ in general relativity
\cite{HE}. The flat chart is identified with a Lorentzian manifold
$$
\mathbb{R}^{3}_{1}(c^2):=(\mathbb{R}^{3},-(dx^{0})^2+e^{2cx^0}
\{(dx^{1})^2+(dx^2)^2\})
$$
of constant sectional curvature $c^2$. This expression shows that
the flat chart is a warped product $\mathbb{E}^{1}_{1}\times_{f}
\mathbb{E}^{2}$ with warping function $f(x^0)=e^{cx^0}$. In
particular, $\mathbb{S}^{3}_{1}(c^2)_+$ is a Robertson-Walker
spacetime.
\end{Remark}

\begin{Example}(Direct product $\mathbb{E}^{1}\times
{\mathbb R}^2_1(c^2)$) Take $(\mu_1,\mu_2)=(0,c)$ with $c\not=0$.
Then the resulting homogeneous spacetime is $\mathbb{R}^3$ with
metric:
$$
-(dx^{0})^{2}+(dx^{1})^{2}+e^{-2cx^0}(dx^{2})^{2},
$$
or equivalently,
$$
(dx^{1})^{2}-(dx^{0})^{2}+e^{-2cx^0}(dx^{2})^{2},
$$
Hence $G(0,c)$ is identified with the direct product of the real
line $\mathbb{E}^{1}(x^1)$ and the warped product model
$$
{\mathbb R}^2_1(c^2)=(\mathbb{R}^2(x^0,x^2),
-(dx^{0})^{2}+e^{-2cx^{0}}(dx^{2})^{2})
$$
of ${\mathbb S}^{2}_1(c^2)_+$. Here, ${\mathbb R}^2_1(c^2)$ denotes
the flat chart model of ${\mathbb S}^2_1(c^2)$. Thus $G(0,c)$ is
identified with $\mathbb{E}^1 \times {\mathbb R}^{2}_1(c^2)$. Note
that $G(0,c)$ is a warped product with trivial warping function.
\end{Example}

\begin{Example}(Homogeneous spacetime $G(c,-c)$)
Let $\mu_1=c, \mu_2=-c$ with $c\ne 0$. Then the resulting
homogeneous spacetime $G(c,-c)$ is the Minkowski motion group
$E(1,1)$ with the Lorentzian metric:
$$-(dx^{0})^{2}+e^{-2cx^0}(dx^{1})^{2}+e^{2cx^0}(dx^{2})^{2}.$$
\end{Example}
\section{Integral representation formula}

Let $\mathfrak{D}(u,v)$ be a simply connected domain and $\varphi:\mathfrak{D}\longrightarrow G(\mu_1,\mu_2)$ an immersion. Let us write $\varphi(u,v)=(x^0(u,v),x^1(u,v),x^2(u,v))$. Then
\begin{equation}
\begin{aligned}
 \alpha'&=\varphi^{-1}\frac{\partial\varphi}{\partial u}\\
 &=\frac{\partial x^0}{\partial u}E_0+\frac{\partial x^1}{\partial u}e^{-\mu_1x^0}E_1+\frac{\partial x^2}{\partial u}e^{-\mu_2x^0}E_2
 \end{aligned}
\end{equation}
and
\begin{equation}
\begin{aligned}
 \alpha''&=\varphi^{-1}\frac{\partial\varphi}{\partial v}\\
 &=\frac{\partial x^0}{\partial v}E_0+\frac{\partial x^1}{\partial v}e^{-\mu_1x^0}E_1+\frac{\partial x^2}{\partial v}e^{-\mu_2x^0}E_2.
 \end{aligned}
\end{equation}
It follows from \eqref{eq:harm3} that
\begin{Lemma}
\label{lem:harm}
 $\varphi$ is harmonic if and only if it satisfies the following equations:
 \begin{equation}
  \begin{aligned}
   \frac{\partial^2x^0}{\partial u\partial v}-\left(\mu_1\frac{\partial x^1}{\partial u}\frac{\partial x^1}{\partial v}e^{-2\mu_1x^0}+\mu_2\frac{\partial x^2}{\partial u}\frac{\partial x^2}{\partial v}e^{-2\mu_2x^0}\right)&=0,\\
   \frac{\partial^2x^1}{\partial u\partial v}-\mu_1\left(\frac{\partial x^0}{\partial u}\frac{\partial x^1}{\partial v}+\frac{\partial x^0}{\partial v}\frac{\partial x^1}{\partial u}\right)&=0,\\
   \frac{\partial^2x^2}{\partial u\partial v}-\mu_2\left(\frac{\partial x^0}{\partial u}\frac{\partial x^2}{\partial v}+\frac{\partial x^0}{\partial v}\frac{\partial x^2}{\partial u}\right)&=0.
  \end{aligned}
 \end{equation}
\end{Lemma}
The exterior derivative $d$ is decomposed as
$$d=\partial'+\partial'',$$
where $\partial'=\frac{\partial}{\partial u}du$ and $\partial''=\frac{\partial}{\partial v}dv$ with respect to the conformal structure of $\mathfrak{D}$. Let
\begin{align*}
 (\omega^0)'&=\frac{\partial x^0}{\partial u}du=\partial' x^0,\\ (\omega^0)''&=\frac{\partial x^0}{\partial v}dv=\partial'' x^0,\\
 (\omega^1)'&=e^{-\mu_1x^0}\partial' x^1,\ (\omega^2)'=e^{-\mu_2x^0}\partial' x^2,\\
 (\omega^1)''&=e^{-\mu_1x^0}\partial'' x^1,\ (\omega^2)''=e^{-\mu_2x^0}\partial'' x^2.
\end{align*}
Then by Lemma \ref{lem:harm}, the 1-forms $(\omega_i)'$, $(\omega_i)''$, $i=0,1,2$ satisfy the differential system:
\begin{align}
\label{eq:harm4}
 \partial''(\omega^0)'&=\mu_1(\omega^1)''\wedge(\omega^1)'+\mu_2(\omega^2)''\wedge(\omega^2)',\\
 \label{eq:harm5}
 \partial''(\omega^i)'&=\mu_i(\omega^i)''\wedge(\omega^0)',\ i=1,2,\\
 \label{eq:harm6}
 \partial'(\omega^0)''&=\mu_1(\omega^1)'\wedge(\omega^1)''+\mu_2(\omega^2)'\wedge(\omega^2)'',\\
 \label{eq:harm7}
 \partial'(\omega^i)''&=\mu_i(\omega^i)'\wedge(\omega^0)'',\ i=1,2.
\end{align}
\begin{Proposition}
If $(\omega_i)'$, $(\omega_i)''$, $i=0,1,2$ satisfy \eqref{eq:harm4}-\eqref{eq:harm7} on a simply connected domain $\mathfrak{D}$. Then
\begin{equation}
\label{eq:integral}
 \varphi(u,v)=\int((\omega^0)',e^{\mu_1x^0}(\omega^1)',e^{\mu_2x^0}(\omega^2)')+\int((\omega^0)'',e^{\mu_1x^0}(\omega^1)'',e^{\mu_2x^0}(\omega^2)'')
\end{equation}
is a harmonic map into $G(\mu_1,\mu_2)$.
\end{Proposition}
Conversely, if $\{(\omega_i)',(\omega_i)'': i=0,1,2\}$ is a solution to \eqref{eq:harm4}-\eqref{eq:harm7} and
\begin{equation}
\label{eq:conformal}
 \begin{aligned}
  -(\omega^0)'\otimes(\omega^0)'+(\omega^1)'\otimes(\omega^1)'+(\omega^2)'\otimes(\omega^2)'&=0,\\
  -(\omega^0)''\otimes(\omega^0)''+(\omega^1)''\otimes(\omega^1)''+(\omega^2)''\otimes(\omega^2)''&=0
 \end{aligned}
\end{equation}
on a simply connected domain $\mathfrak{D}$, then $\varphi(u,v)$ in \eqref{eq:integral} is a weakly conformal harmonic map into $G(\mu_1,\mu_2)$. In addition, if 
\begin{equation}
 -(\omega^0)'\otimes(\omega^0)''+(\omega^1)'\otimes(\omega^1)''+(\omega^2)'\otimes(\omega^2)''\ne 0,
\end{equation}
then $\varphi(u,v)$ in \eqref{eq:integral} is a minimal timelike surface in $G(\mu_1,\mu_2)$.
\section{Normal Gau{\ss} map}
Let $\varphi:\mathfrak{D}\longrightarrow G(\mu_1,\mu_2)$ be a Lorentz surface i.e. a conformal timelike surface. Take a unit normal vector field $N$ along $\varphi$. Then by the left translation we obtain the smooth map
$$\varphi^{-1}N:\mathfrak{D}\longrightarrow\mathbb{S}^2_1(1),$$
where
$$\mathbb{S}^2_1(1)=\{u^0E_0+u^1E_1+u^2E_2: -(u^0)^2+(u^1)^2+(u^2)^2=1\}\subset\mathfrak{g}(\mu_1,\mu_2)$$
is the de Sitter 2-space of constant Gau\ss ian curvature 1. The Lie algebra $\mathfrak{g}(\mu_1,\mu_2)$ is identified with Minkowski 3-space $\mathbb{E}^3_1(u^0,u^1,u^2)$ via the orthonormal basis $\{E_0,E_1,E_2\}$. Then smooth map $\varphi^{-1}N$ is called the normal Gau{\ss} map of $\varphi$. Let $\varphi:\mathfrak{D}\longrightarrow G(\mu_1,\mu_2)$ be a minimal timelike surface determined by the data $((\omega^0)',(\omega^1)',(\omega^2)')$ and $((\omega^0)'',(\omega^1)'',(\omega^2)'')$. Write $(\omega^i)'=\xi^idu$ and $(\omega^i)''=\eta^idv$, $i=0,1,2$. Then
\begin{equation}
 \begin{aligned}
  I&=2(-(\omega^0)'\otimes(\omega^0)''+(\omega^1)'\otimes(\omega^1)''+(\omega^2)'\otimes(\omega^2)'')\\
  &=2(-\xi^0\eta^0+\xi^1\eta^1+\xi^2\eta^2)dudv.
 \end{aligned}
\end{equation}
The conformality condition \eqref{eq:conformal} can be written as
\begin{equation}
\label{eq:conformal2}
 \begin{aligned}
  -(\xi^0)^2+(\xi^1)^2+(\xi^2)^2&=0,\\
  -(\eta^0)^2+(\eta^1)^2+(\eta^2)^2&=0.
 \end{aligned}
\end{equation}
It follows from \eqref{eq:conformal2} that one can introduce pairs of functions $(q,f)$ and $(r,g)$ such that
\begin{equation}
 \begin{aligned}
  q&=\frac{-\xi^2}{\xi^0-\xi^1},\ f=\xi^0-\xi^1,\\
  r&=\frac{\eta^2}{\eta^0+\eta^1},\ g=-(\eta^0+\eta^1).
 \end{aligned}
\end{equation}
In terms of $(q,f)$ and $(r,g)$, $\varphi(u,v)=(x^0(u,v),x^1(u,v),x^2(u,v))$ is given by Weierstra{\ss} type representation formula
\begin{equation}
\label{eq:weierstrass}
 \begin{aligned}
  x^0(u,v)&=\frac{1}{2}\int(1+q^2)fdu-(1+r^2)gdv,\\
  x^1(u,v)&=-\frac{1}{2}e^{\mu_1x^0(u,v)}\int(1-q^2)fdu+(1-r^2)gdv,\\
  x^2(u,v)&=-e^{\mu_2x^0(u,v)}\int qfdu+rgdv.
 \end{aligned}
\end{equation}
with first fundamental form
\begin{equation}
 I=(1+qr)^2fgdudv.
\end{equation}
\begin{Remark}
In the study of minimal timelike surfaces in Minkowski 3-space, one may assume that $f=g=1$ so that \eqref{eq:weierstrass} reduces to a simpler form called the \emph{normalized Weierstra{\ss} formula}. This is possible as there are no restrictions on $f$ and $g$ other than $f$ and $g$ being Lorentz holomorphic and Lorentz anti-holomorphic respectively. (See \cite{Inoguchi-Toda} and \cite{Lee2}.) However, this is not the case with minimal timelike surfaces in de Sitter 3-space as we will see later.
\end{Remark}
It turns out that the pair $(q,r)$ is the Normal Gau{\ss} map $\varphi^{-1}N$ projected into the Minkowski 2-pane $\mathbb{E}^2_1$. To see this, first the normal Gau{\ss} map is computed to be
\begin{equation}
 \varphi^{-1}N=\frac{1}{qr+1}[(q-r)E_0+(q+r)E_1+(qr-1)E_2].
\end{equation}
Let $\wp_\mathcal{N}:\mathbb{S}^2_1(1)\setminus\{x^2=1\}\longrightarrow\mathbb{E}^2_1\setminus\mathbb{H}^1_0$ be the stereographic projection from the north pole $\mathcal{N}=(0,0,1)$. Here, $\mathbb{H}^1_0$ is the hyperbola
$$\mathbb{H}^1_0=\{x^0E_0+x^1E_1\in\mathbb{E}^2_1:-(x^0)^2+(x^1)^2=-1\}.$$
Then
\begin{equation}
\wp_\mathcal{N}(x^0E_0+x^1E_1+x^2E_2)=\frac{x^0}{1-x^2}E^0+\frac{x^1}{1-x^2}E^1.
\end{equation}
So, the normal Gau{\ss} map $\varphi^{-1}N$ is projected into the Minkowski plane $\mathbb{E}^2_1$ via $\wp_\mathcal{N}$ as
\begin{equation}
\label{eq:projectgauss}
 \wp_\mathcal{N}\circ\varphi^{-1}N=\frac{q-r}{2}E_0+\frac{q+r}{2}E_1\in\mathbb{E}^2_1(t,x).
\end{equation}
In terms of null coordinates $(u,v)$, \eqref{eq:projectgauss} is written as
\begin{equation}
 \wp_\mathcal{N}\circ\varphi^{-1}N=(q,r)\in\mathbb{E}^2_1(u,v).
\end{equation}
The pair $(q,r)$ is called the \emph{projected normal Gau{\ss} map} of $\varphi$.
It follows from \eqref{eq:harm4} and \eqref{eq:harm5} that
\begin{equation}
\label{eq:harm8}
 \begin{aligned}
  \frac{\partial\xi^0}{\partial v}&=\mu_1\eta^1\xi^1+\mu_2\eta^2\xi^2,\\
  \frac{\partial\xi^i}{\partial v}&=\mu_i\eta^i\xi^0,\ i=1,2.
 \end{aligned}
\end{equation}
Using \eqref{eq:harm8}, we obtain
\begin{equation}
\label{eq:projectgauss2a}
\begin{aligned}
  \frac{\partial f}{\partial v}&=\frac{\partial\xi^0}{\partial v}-\frac{\partial\xi^1}{\partial v}\\
  &=\frac{\mu_1}{2}(1-r^2)fg+\mu_2qrfg
  \end{aligned}
  \end{equation}
  and
  \begin{equation}
  \label{eq:projectgauss2b}
  \begin{aligned}
  \frac{\partial q}{\partial v}&=-\frac{\frac{\partial\xi^2}{\partial v}f-\xi^2\frac{\partial f}{\partial v}}{f^2}\\
  &=-\frac{\mu_1}{2}q(1-r^2)g+\frac{\mu_2}{2}(1-q^2)rg.
\end{aligned}
  \end{equation}
It follows from \eqref{eq:harm6} and \eqref{eq:harm7} that
\begin{equation}
\label{eq:harm9}
 \begin{aligned}
  \frac{\partial\eta^0}{\partial u}&=\mu_1\xi^1\eta^1+\mu_2\xi^2\eta^2,\\
  \frac{\partial\eta^i}{\partial u}&=\mu_i\xi^i\eta^0,\ i=1,2.
 \end{aligned}
\end{equation}
Using \eqref{eq:harm9}, we obtain
\begin{equation}
\begin{aligned}
\label{eq:projectgauss3a}
  \frac{\partial g}{\partial u}&=-\frac{\partial\eta^0}{\partial u}-\frac{\partial\eta_1}{\partial u}\\
  &=-\frac{\mu_1}{2}(1-q^2)fg-\mu_2qrfg
  \end{aligned}
  \end{equation}
  and
  \begin{equation}
\begin{aligned}
  \label{eq:projectgauss3b}
  \frac{\partial r}{\partial u}&=-\frac{\frac{\partial\eta^2}{\partial u}g-\eta^2\frac{\partial g}{\partial u}}{g^2}\\
  &=\frac{\mu_1}{2}(1-q^2)rf-\frac{\mu_2}{2}q(1-r^2)f.
 \end{aligned}
 \end{equation}
\begin{Remark}
\label{rem:desitter2}
 Setting $f=g=1$, we obtain from \eqref{eq:projectgauss2a}, \eqref{eq:projectgauss2b}, \eqref{eq:projectgauss3a}, and \eqref{eq:projectgauss3b}
\begin{align}
\label{eq:projectgauss4a}
\mu_2qr&=-\frac{\mu_1}{2}(1-r^2),\\
\label{eq:projectgauss4b}
  \frac{\partial q}{\partial v}&=-\frac{\mu_1}{2}(1-r^2)q+\frac{\mu_2}{2}(1-q^2)r
\end{align}
and
\begin{align}
 \label{eq:projectgauss5a}
\mu_2qr&=-\frac{\mu_1}{2}(1-q^2),\\
\label{eq:projectgauss5b}
  \frac{\partial r}{\partial u}&=\frac{\mu_1}{2}(1-q^2)r-\frac{\mu_2}{2}q(1-r^2).
\end{align}
It follows from \eqref{eq:projectgauss4a} and \eqref{eq:projectgauss5a} that $q=\pm r$. Let $\mu_1=\mu_2=c\ne 0$. If $q=r$ then $\frac{\partial q}{\partial v}=\frac{\partial r}{\partial u}=0$. This means that $q=r$ is a constant, say $A$. By \eqref{eq:weierstrass} $\varphi$ is computed to be
$$
\begin{aligned}
\varphi(u,v)=&\left(\frac{1}{2}(1+A^2)(u-v),-\frac{1}{2}e^{\frac{1}{2}c(1+A^2)(u-v)}(1-A^2)(u+v),\right.\\
&\left.-e^{\frac{1}{2}c(1+A^2)(u-v)}(u+v)\right)
\end{aligned}$$
or
$$\varphi(t,x)=((1+A^2)t,-e^{c(1+A^2)t}x,-2e^{c(1+A^2)t}Ax).$$
This surface cannot be minimal as it is not conformal. If $q=-r\ne 0$ then from \eqref{eq:projectgauss4b} and \eqref{eq:projectgauss5b} we obtain the separable differential equations
\begin{align}
\label{eq:de1}
 \frac{1}{q(1-q^2)}\frac{\partial q}{\partial v}&=-c,\\
 \label{eq:de2}
 \frac{1}{r(1-r^2)}\frac{\partial r}{\partial u}&=c.
\end{align}
\eqref{eq:de1} has solution
\begin{equation}
\label{eq:sol1}
 q\sqrt{\frac{1-q}{1+q}}=A(u)e^{-cv},
\end{equation}
where $A(u)>0$ is a Lorentz holomorphic function. \eqref{eq:de2} has solution
\begin{equation}
\label{eq:sol2a}
 r\sqrt{\frac{1-r}{1+r}}=B(v)e^{cu},
\end{equation}
where $B(v)>0$ is a Lorentz anti-holomorphic function. Since $q=-r$, \eqref{eq:sol2a} can be written as
\begin{equation}
 \label{eq:sol2b}
 -q\sqrt{\frac{1+q}{1-q}}=B(v)e^{cu}.
\end{equation}
\eqref{eq:sol1} and \eqref{eq:sol2b} yield
$$q^2=-A(u)B(v)e^{c(u-v)}<0.$$
This case cannot occur as $q$ is a real-valued function.
\end{Remark}
As seen in Section \ref{sec:sol}, $G(0,0)=\mathbb{E}^3_1$ and $G(c,c)=\mathbb{S}^3_1(c^2)_+$ are the only cases of solvable Lie group $G(\mu_1,\mu_2)$ with constant sectional curvature.
\begin{Remark}
For $G(0,0)=\mathbb{E}^3_1$,
\begin{align*}
 \frac{\partial f}{\partial v}&=\frac{\partial q}{\partial v}=0,\\
 \frac{\partial g}{\partial u}&=\frac{\partial r}{\partial u}=0.
\end{align*}
That is, $f,q$ are Lorentz holomorphic and $g,r$ are Lorentz anti-holomorphic. From \eqref{eq:weierstrass}, we retrieve the Weierstra{\ss} representation formula (\cite{Inoguchi-Toda}, \cite{Lee2}) for minimal timelike surface $\varphi(u,v)=(x^0(u,v),x^1(u,v),x^2(u,v))$ in $\mathbb{E}^3_1$ given by
\begin{equation}
\label{eq:weierstrass2}
 \begin{aligned}
  x^0(u,v)&=\frac{1}{2}\int(1+q^2)fdu-(1+r^2)gdv,\\
  x^1(u,v)&=-\frac{1}{2}\int(1-q^2)fdu+(1-r^2)gdv,\\
  x^2(u,v)&=-\int qfdu+rgdv.
 \end{aligned}
\end{equation}
\end{Remark}
\begin{Remark}
If $\mu_1=\mu_2=c\ne 0$, then \eqref{eq:projectgauss2b} and \eqref{eq:projectgauss3b} can be written respectively as
\begin{align}
 \frac{\partial q}{\partial v}&=\frac{c}{2}g(r-q)(1+qr),\\
 \frac{\partial r}{\partial u}&=\frac{c}{2}f(r-q)(1+qr).
\end{align}
If $q$ is Lorentz holomorphic, then $q=r$ or $1+qr=0$. If $1+qr=0$ then $I=0$. $q=r$ cannot occur as discussed in Remark \ref{rem:desitter2}. Hence, $q$ cannot be Lorentz holomorphic for minimal timelike surfaces in $\mathbb{S}^3_1(c^2)_+$. For the same reason, $r$ cannot be Lorentz anti-holomorphic for minimal timelike surfaces in $\mathbb{S}^3_1(c^2)_+$.
\end{Remark}
From here on, we assume that $q^2\ne 1$ and $r^2\ne 1$. It follows from \eqref{eq:projectgauss2a}, \eqref{eq:projectgauss2b}, \eqref{eq:projectgauss3a}, and \eqref{eq:projectgauss3b} that the projected normal Gau{\ss} map $(q,r)$ satisfies the equations
\begin{equation}
\label{eq:harmeq1a}
\begin{aligned}
 &\frac{\partial^2q}{\partial u\partial v}+\frac{\mu_1(1-r^2)+2\mu_2qr}{-\mu_1q(1-r^2)+\mu_2(1-q^2)r}\frac{\partial q}{\partial u}\frac{\partial q}{\partial v}\\
 &+\frac{(\mu_1^2-\mu_2^2)(1-q^2)(1+r^2)q}{[-\mu_1q(1-r^2)+\mu_2(1-q^2)r][-\mu_1(1-q^2)r+\mu_2q(1-r^2)]}\frac{\partial r}{\partial u}\frac{\partial q}{\partial v}\\
 &=0
\end{aligned}
\end{equation}
and
\begin{equation}
\label{eq:harmeq1b}
 \begin{aligned}
  &\frac{\partial^2r}{\partial v\partial u}+\frac{\mu_1(1-q^2)+2\mu_2qr}{-\mu_1(1-q^2)r+\mu_2q(1-r^2)}\frac{\partial r}{\partial u}\frac{\partial r}{\partial v}\\
 &+\frac{(\mu_1^2-\mu_2^2)(1+q^2)(1-r^2)r}{[-\mu_1(1-q^2)r+\mu_2q(1-r^2)][-\mu_1q(1-r^2)+\mu_2(1-q^2)r]}\frac{\partial r}{\partial u}\frac{\partial q}{\partial v}\\
 &=0.
 \end{aligned}
\end{equation}
The equations \eqref{eq:harmeq1a} and \eqref{eq:harmeq1b} are not the harmonic map equations for the projected normal Gau{\ss} map $(q,r)$ in general. The following theorem tells under what conditions they become the harmonic map equations for $(q,r)$.
\begin{Theorem}
\label{thm:harmeq}
 The projected normal Gau{\ss} map $(q,r)$ is a harmonic map if and only if $\mu_1^2=\mu_2^2$. If $\mu_1=\mu_2\ne 0$ then \eqref{eq:harmeq1a} and \eqref{eq:harmeq1b} reduce to
 \begin{align}
 \label{eq:harmeq2a}
  \frac{\partial^2q}{\partial u\partial v}+\frac{1-r^2+2qr}{(1-q^2)r-q(1-r^2)}\frac{\partial q}{\partial u}\frac{\partial q}{\partial v}&=0,\\
\label{eq:harmeq2b}
  \frac{\partial^2 r}{\partial v\partial u}+\frac{-(1-q^2)-2qr}{(1-q^2)r-q(1-r^2)}\frac{\partial r}{\partial u}\frac{\partial r}{\partial v}&=0.
 \end{align}
 \eqref{eq:harmeq2a} and \eqref{eq:harmeq2b} are the harmonic map equations for the map $(q,r):\mathfrak{D}(u,v)\longrightarrow
\left(\mathbb{E}^2_1(\alpha,\beta),\frac{2d\alpha d\beta}{(1-\alpha^2)\beta-\alpha(1-\beta^2)}\right)$. If $\mu_1=-\mu_2$ then \eqref{eq:harmeq1a} and \eqref{eq:harmeq1b} reduce to
 \begin{align}
 \label{eq:harmeq3a}
  \frac{\partial^2q}{\partial u\partial v}+\frac{-(1-r^2)+2qr}{(1-q^2)r+q(1-r^2)}\frac{\partial q}{\partial u}\frac{\partial q}{\partial v}&=0,\\
\label{eq:harmeq3b}
  \frac{\partial^2 r}{\partial v\partial u}+\frac{-(1-q^2)+2qr}{(1-q^2)r+q(1-r^2)}\frac{\partial r}{\partial u}\frac{\partial r}{\partial v}&=0.
 \end{align}
 \eqref{eq:harmeq3a} and \eqref{eq:harmeq3b} are the harmonic map equations for the map $(q,r):\mathfrak{D}(u,v)\longrightarrow
\left(\mathbb{E}^2_1(\alpha,\beta),\frac{2d\alpha d\beta}{(1-\alpha^2)\beta+\alpha(1-\beta^2)}\right)$.
\end{Theorem}
{\it Proof.} The tension field $\tau(q,r)$ of $(q,r)$ is given by (\cite{EL}, \cite{Wood})
\begin{equation}
 \tau(q,r)=4\lambda^{-2}\left(\frac{\partial^2q}{\partial u\partial v}+\Gamma^\alpha_{\alpha\alpha}\frac{\partial q}{\partial u}\frac{\partial q}{\partial v},\frac{\partial^2r}{\partial v\partial u}+\Gamma^\beta_{\beta\beta}\frac{\partial r}{\partial u}\frac{\partial r}{\partial v}\right),
\end{equation}
where $\lambda$ is a parameter of conformality. Here, $\Gamma^\alpha_{\alpha\alpha},\Gamma^\beta_{\beta\beta}$ are the Christoffel symbols of $\mathbb{E}^2_1(\alpha,\beta)$. Comparing \eqref{eq:harmeq1a}, \eqref{eq:harmeq1b} and $\tau(q,r)=0$, we see that \eqref{eq:harmeq1a} and \eqref{eq:harmeq1b} are the harmonic map equations for $(q,r)$ if and only if $\mu_1^2=\mu_2^2$. In order to find a metric on $\mathbb{E}^2_1(\alpha,\beta)$ with which \eqref{eq:harmeq1a} and \eqref{eq:harmeq1b} are the harmonic map equations, one needs to solve the first-order partial differential equations
\begin{equation}
\label{eq:metric}
 \begin{aligned}
  \Gamma^\alpha_{\alpha\alpha}&=g^{\alpha\beta}\frac{\partial g_{\alpha\beta}}{\partial\alpha}\\
  &=\left\{\begin{array}{ccc}
               \frac{1-\beta^2+2\alpha\beta}{(1-\alpha^2)\beta-\alpha(1-\beta^2)} & \mbox{if} & \mu_1=\mu_2\ne 0,\\
               \\
               \frac{-(1-\beta^2)+2\alpha\beta}{(1-\alpha^2)\beta+\alpha(1-\beta^2)} & \mbox{if} & \mu_1=-\mu_2,
               \end{array}\right.\\
   \Gamma^\beta_{\beta\beta}&=g^{\alpha\beta}\frac{\partial g_{\alpha\beta}}{\partial\beta}\\ 
   &=\left\{\begin{array}{ccc}
               \frac{-(1-\alpha^2)-2\alpha\beta}{(1-\alpha^2)\beta-\alpha(1-\beta^2)} & \mbox{if} & \mu_1=\mu_2\ne 0,\\
               \\
               \frac{-(1-\alpha^2)+2\alpha\beta}{(1-\alpha^2)\beta+\alpha(1-\beta^2)} & \mbox{if} & \mu_1=-\mu_2.
               \end{array}\right.
 \end{aligned}
\end{equation}
The solutions are given by
\begin{equation}
(g_{\alpha\beta})=\left\{\begin{array}{ccc}
                    \begin{pmatrix}
                     0 & \frac{1}{(1-\alpha^2)\beta-\alpha(1-\beta^2)}\\
                     \frac{1}{(1-\alpha^2)\beta-\alpha(1-\beta^2)} & 0
                    \end{pmatrix} & \mbox{if} & \mu_1=\mu_2\ne 0,\\
                    \\
                    \begin{pmatrix}
                     0 & \frac{1}{(1-\alpha^2)\beta+\alpha(1-\beta^2)}\\
                     \frac{1}{(1-\alpha^2)\beta+\alpha(1-\beta^2)} & 0
                    \end{pmatrix} & \mbox{if} & \mu_1=-\mu_2.
                   \end{array}\right.
\end{equation}
Q.E.D.
\begin{Remark}
 Clearly, the projected normal Gau{\ss} map $(q,r)$ of a minimal timelike surface in $G(0,0)=\mathbb{E}^3_1$ satisfies the wave equation
 \begin{equation}
  \Box(q,r)=0,
 \end{equation}
 where $\Box$ denotes the d'Alembertian
 \begin{equation}
  \Box=\lambda^{-2}\left(-\frac{\partial^2}{\partial t^2}+\frac{\partial^2}{\partial x^2}\right)=4\lambda^{-2}\frac{\partial^2}{\partial u\partial v}.
 \end{equation}
\end{Remark}
\begin{Remark}
 Theorem \ref{thm:harmeq} tells that Minkowski 3-space $G(0,0)=\mathbb{E}^3_1$, de Sitter 3-space $G(c,c)=\mathbb{S}^3_1$, and $G(c,-c)=E(1,1)$ are the only homogeneous 3-dimensional spacetimes among $G(\mu_1,\mu_2)$ in which the projected normal Gau{\ss} map of a minimal timelike surface is harmonic.
\end{Remark}

\noindent {\sc Sungwook Lee\\ 
Department of Mathematics\\
University of Southern Mississippi\\
118 College Drive, \#5045\\
Hattiesburg, MS 39406-0001, U.S.A.}

\smallskip

\noindent {\it E-mail address}: {\tt sunglee@usm.edu}


\begin{thebibliography}{9999}
\bibitem{EL}
J. Eells and L. Lemaire, Selected topics in harmonic maps, {\it C.B.M.S. Regional Conference Series} {\bf 50}, Amer. Math. Soc. 1983.

\bibitem{HE}
S.~W.~Hawking and G.~F.~R.~Ellis,
{\it The Large Scale Structure of Space-Time},
Cambridge Univ. Press, Cambridge, 1973.

\bibitem{Inoguchi-Toda}
J. Inoguchi and M. Toda, Timelike minimal surfaces via loop groups, {\it Acta. Appl. Math.} {\bf 83} (2004), 313--335

\bibitem{Kobayashi}
O.~Kobayashi, Maximal surfaces in the 3-dimensional Minkowski space,
{\it Tokyo J. Math.} {\bf 6} (1983), 297--309.

\bibitem{Lee}
S. Lee, Maximal surfaces in a certain 3-dimensional homogeneous spacetime, {\it Differential Geometry and Its Applications} {\bf 26} (2008), Issue 5, 536--543.

\bibitem{Lee2}
S. Lee, Weierstrass representation for timelike minimal surfaces in Minkowski 3-space, {\it Communications in Mathematical Analysis}, {\bf Conf. 01}(2008), 11--19

\bibitem{Mc}
L.~McNertney,
{\sl
One-parameter families of surfaces with constant
curvature in Lorentz 3-space},
Ph.~D. Thesis, Brown Univ., Providence, RI, U.S.A., 1980.

\bibitem{Thurston}
W. M. Thurston, Three-dimensional Geometry and Topology I, Princeton Math. Series., vol. {\bf 35} (S. Levy ed.), 1997.

\bibitem{Wood}
J. C. Wood, Harmonic maps into symmetric spaces and integrable systems, Aspects of Mathematics, vol. {\bf E23}, Vieweg, Braunschweig/Wiesbaden, 1994, 29--55.
\end{thebibliography}
\end{document}